\documentstyle[12pt]{article}
\input{amssym.def}
\input{amssym.tex}
\input{cyracc.def}
\newfont{\cyrfnt}{wncyr10}
\newcommand{\cyr}{\baselineskip12.5pt\cyrfnt\cyracc}
\newfont{\cyifnt}{wncyi10}
\newcommand{\cyi}{\baselineskip12.5pt\cyifnt\cyracc}
\newtheorem{thm}{}[section]
\author{J\"org Jahnel\footnotemark}
\date{Revised version, August 1995}
\title{Heights for line bundles on arithmetic surfaces}

\begin{document}
\renewcommand{\thefootnote}{\fnsymbol{footnote}}
\maketitle
\begin{abstract}
For line bundles on arithmetic varieties we construct height functions using
arithmetic intersection theory. In the case of an arithmetic surface,
generically of genus
$g$,
for line bundles of degree
$g$
equivalence is shown to the height on the Jacobian defined by
$\Theta$.
\end{abstract}

\section{Introduction}
\footnotetext[1]{Mathematisches Institut, Universit\"at G\"ottingen,
Bunsenstra{\ss}e 3-5, 37073 G\"ottingen, Germany
\newline
email: jahnel@cfgauss.uni-math.gwdg.de}
\thispagestyle{empty}
In this paper we will suggest a construction for height functions for line
bundles on arithmetic varieties. Following the philosophy of
\cite{Bost/Gillet/Soule 93} heights should be objects in arithmetic geometry
analogous to degrees in algebraic geometry. So let
$K$
be a number field,
${\cal O}_{K}$
its ring of integers and
${\cal X} / {\cal O}_{K}$
an arithmetic variety, i.e. a regular scheme, projective and flat over
${\cal O}_{K}$,
whose generic fiber
$X / K$
we assume to be connected of dimension
$d$.
Then we have to fix a metrized line bundle
$( {\cal T}, \| . \| )$
or, equivalently, its first Chern class
$$\stackrel{\wedge}{{c}_{1}} ({\cal T} , \| . \| ) = (T,g_{T}) \in ~
\stackrel{\wedge}{{ \rm CH}^{1}} ({\cal X}) ~~.$$

The height of a line bundle
${\cal L}$
on
${\cal X}$
should be the arithmetic degree of the intersection of
$\stackrel{\wedge}{{c}_{1}}({\cal L} )$
with
$(T,g_{T})^{d}$.
For this a natural hermitian metric has to be chosen on
${\cal L}$.
We fix a K\"ahler metric
$\omega_{0}$
on
${\cal X} ({\Bbb{C}})$,
invariant under complex conjugation
$F_{\infty}$,
as in \cite{Arakelov 74}. Then it is well known that the condition on the Chern
form to be harmonic defines
$\| .\|$
up to a locally constant factor.

In order to determine this factor we require
$$\stackrel{\wedge}{\deg} \Big( \det R \pi_{*} {\cal L}, \| . \|_{Q} \Big) = 0
.$$
Here
$\pi : {\cal X} \longrightarrow {\rm Spec}~ {\cal O}_{K}$
is the structural morphism and
$\| . \|_{Q}$
is Quillen's metric (\cite{Quillen 85}, \cite{Bismut/Gillet/Soule 88}) at the
infinite places of
$K$.

\begin{thm} {\bf Fact.} {\rm a)} If the Euler characteristic
$\chi ({\cal L})$
does not vanish, such a metric exists.
\newline
{\rm b)}
$\stackrel{\wedge}{{c}_{1}} ({\cal L}, \| .\|)$
is uniquely determined up to a summand
$(0,C)$,
where
$C = (C_{\sigma})_{\sigma : K \hookrightarrow \Bbb{C}}$
is a system of constants on
$X \times_{{\rm Spec K}, \sigma} {\rm Spec} ~ \Bbb{C}$
with
$$\sum_{\sigma : K \hookrightarrow \Bbb{C}} C_{\sigma} = 0 ~~~~~(and
{}~~C_{\sigma} = C_{\bar{\sigma}}).$$
\end{thm}

\begin{thm} {\bf Fact.}
Such
$(0,C) \in ~ \stackrel{\wedge}{{\rm CH}^{1}} ({\cal X})$
are numerically trivial.
\end{thm}

\begin{thm}
{\rm Now we can state our fundamental} \newline
{\bf Definition.}
The {\rm height} of the line bundle
${\cal L}$
is given by
$$h_{{\cal T}, \omega_{0}} ({\cal L}) := ~ \stackrel{\wedge}{\deg} ~\pi_{*} \,
\Big[ \stackrel{\wedge}{{c}_{1}} ({\cal L}, \| . \|) \cdot (T, g_{T})^{d} \Big]
,$$
where
$\| . \|$
is one of the distinguished metrics specified above.
\end{thm}

\begin{thm}\label{main}
{\rm In this paper we will analyze this definition in the case of arithmetic
surfaces. Our main result is} \newline
{\bf Theorem.}
Let
${\cal C} / {\cal O}_{K}$
be a regular projective variety of dimension
$2$,
flat over
${\cal O}_{K}$
and generically connected of genus
$g$,
$x \in ({\cal C} \times_{{\rm Spec} {\cal O}_{K}} {\rm Spec} ~ K) (K)$
be a
$K$-valued
point and
$\Theta$
be the Theta divisor on the Jacobian
$J = {\rm Pic}^{g} (C)$
(defined using
$x$).
On
$$\coprod_{\sigma: K \hookrightarrow \Bbb{C}} \Big( {\cal C} \times_{{\rm Spec}
{\cal O}_{K}, \sigma} {\rm Spec} ~ \Bbb{C} \Big) (\Bbb{C})$$
let
$\omega$
be a K\"ahler form invariant under
$F_{\infty}$
and normalized by
$$\int_{({\cal C} \times_{{\rm Spec} {\cal O}_{K}, \sigma} {\rm Spec}
\Bbb{C})(\Bbb{C})} \omega = 1$$
for every
$\sigma$.

Then, for line bundles
${\cal L} / {\cal C}$,
fiber-by-fiber of degree
$g$
and of degree of absolute value less than
$H$
on every irreducible component of the special fibers of
${\cal C}$
(with some constant
$H \in \Bbb{N}$)
$$h_{x,\omega} ({\cal L}) = h_{\Theta} ({\cal L}_{K}) + {\rm O}(1) ,$$
where
$h_{\Theta}$
is the height on
$J$
defined using the ample divisor
$\Theta$.
\end{thm}

\begin{thm} {\bf Remark.}
{\rm Another connection between heights on the Jacobian of a curve and
arithmetic intersection theory was obtained by Faltings \cite{Faltings 84} and
Hriljac \cite{Hriljac 85}. Recently it has been generalized to higher
dimensions and higher codimension Chow groups by K\"unnemann \cite{Kunnemann
95}. They can write down an explicit formula for the N\'{e}ron-Tate height
pairing on the Jacobian (higher Picard variety) in terms of arithmetic
intersection theory. The main point is that they consider line bundles (cycles)
algebraically equivalent to zero. So there is no need for them to scale a
metric (to specify the infinite part of the arithmetic cycles occuring). Our
approach, to the contrary, seems to work best for sufficiently ample algebraic
equivalence classes of line bundles. A formal relationship between our approach
and the other one is not known to the author.}
\end{thm}

\begin{thm} {\rm In order to prove the two facts above we will use the
following simple} \newline
{\bf Lemma.}
Let
$f: X \longrightarrow Y$
be a smooth proper map of complex manifolds, where
$X$
has a K\"ahler structure
$\omega$
and
$Y$
is connected, and
$E$
be a holomorphic vector bundle on
$X$.
For a hermitian metric
$\| . \|$
on
$E$
and a constant factor
$D > 0$
we have
$$ h_{Q,(E,D \cdot \| . \|)} = h_{Q,(E,\| . \|)} \cdot D^{\chi (E)} .$$
\end{thm}
{\bf Proof.}
The homomorphism
\begin{eqnarray*}
(E,\| . \|) & \longrightarrow & (E, D \cdot \| . \|) \\
s & \mapsto & \frac{1}{D} \cdot s
\end{eqnarray*}
is an isometry inducing the isometry
\pagebreak
\begin{eqnarray*}
\Big( \det R\pi_{*} E, h_{Q,(E,\| . \|)} \Big) & \longrightarrow & \Big( \det
R\pi_{*} E, h_{Q,(E,D \cdot \| . \|)} \Big) \\
x & \mapsto & D^{-\chi (E)} \cdot x ~~~~~~~~~~~~~~~~~.
\end{eqnarray*}
\begin{center}
$\Box$
\end{center}

\begin{thm}
{\bf Proof of Fact 1. Existence:} {\rm Multiplication of
$\| . \|$
by
$D$
will change the Quillen metric by the factor
$D^{\chi (E)}$
and therefore
$\stackrel{\wedge}{\deg} \det R \pi_{*} E$
by the summand
$[K:{\Bbb{Q}}] \chi(E) \log D$.
\newline
{\bf Uniqueness:} The harmonicity condition and invariance under
$F_{\infty}$
determine
$\| . \|$
up to constant factors
$D_{\sigma} > 0$
for each
$\sigma: K \hookrightarrow \Bbb{C}$
with
$D_{\sigma} = D_{\bar{\sigma}}$.
The scaling condition requires
$$\prod_{\sigma: K \hookrightarrow \Bbb{C}} D_{\sigma}^{\chi(E)} = 1$$
or
$\sum_{\sigma: K \hookrightarrow \Bbb{C}} \log D_{\sigma} = 0.$ }
\begin{center}
$\Box$
\end{center}
\end{thm}

\begin{thm}
{\bf Proof of Fact 2.} {\rm Let
$(Z,g_{Z}) \in ~ \stackrel{\wedge}{\rm CH}_{1}({\cal X})$. Then
\begin{eqnarray*}
(0,C) \cdot (Z,g_{Z}) & = &
(0,g_{Z} \cdot \omega_{(0,C)} + C \cdot \delta_{Z}) \\
  & = & (0,C \cdot \delta_{Z}) ~~~~~~~~.
\end{eqnarray*}
$Z$
is a zero-cycle on
$X$,
so
$\delta_{Z}$
will have, independently on
$\sigma$,
always the integral
$\deg Z$.
Therefore
\begin{eqnarray*}
\stackrel{\wedge}{\deg} ~ \pi_{*} \Big[ (0,C) \cdot (Z,g_{Z}) \Big] & = &
\frac{\scriptstyle 1}{\scriptstyle 2} ~ \sum_{\sigma} ~~ \left[ C_{\sigma}
\int_{X \times_{{\rm Spec K}, \sigma} {\rm Spec} ~ {\Bbb{C}} (\Bbb{C})}
\delta_{Z} \right] \\
 & = & \frac{\scriptstyle 1}{\scriptstyle 2} ~ (\sum_{\sigma} C_{\sigma}) \cdot
\deg Z \\
 & = & 0 ~~.
\end{eqnarray*} }
\begin{center}
$\Box$
\end{center}
\end{thm}

\section{Divisors versus points of the Jacobian}
\begin{thm}
{\rm The remainder of this paper is devoted to the proof of Theorem \ref{main}.
So let
$C/K$
be a regular proper algebraic curve of genus
$g$
with
$C(K) \neq \emptyset$.
We consider a regular projective model
${\cal C} / {\cal O}_{K}$.
Denote by
$J = {\rm Pic}^{g}_{C/K}$
the Jacobian of
$C$.
When
$x \in C(K)$
is chosen we have a canonical isomorphism
${\rm Pic}^{g-1}_{C/K} \longrightarrow {\rm Pic}^{g}_{C/K} = J$
and thus the divisor
$\Theta$
on
$J$.
$\Theta$ induces a closed embedding
$i^{'}: J \hookrightarrow {\bf P}^{N}_{K}$
and a "naive" height for
$K$-valued
points of
$J$:
$$h_{\Theta} (D) := \log ~ \left( \prod_{\nu \in M_{K}} \max \left\{ \|
i(D)_{0} \|_{\nu}, ~ \ldots ~ , \| i(D)_{N} \|_{\nu} \right\} \right) ~~.$$
Accordingly
$j^{*} (\Theta)$
induces a morphism
$i: C^{g} \stackrel{j}{\longrightarrow} J \stackrel{i^{'}}{\longrightarrow}
{\bf P}^{N}_{K}$
and a height function
$h_{j^{*} (\Theta)}$
for
$K$-valued
points of
$C^{g}$.
Here
$j$
denotes the natural map sending a divisor to its associated line bundle. A
general construction for heights defined by a divisor, the "height machine", is
given in [CS, Chapter VI, Theorem 3.3].

The underlying height
$h$
for
$K$-valued points of
${\bf P}^{N}_{K}$
is a height in the sense of Arakelov theory \cite{Bost/Gillet/Soule 93} as
follows: We choose the regular projective model
${\bf P}^{N}_{{\cal O}_{K}} \supseteq {\bf P}^{N}_{K}$.
Every \linebreak
$K$-valued
point
$y$
of
${\bf P}^{N}_{K}$
can be extended uniquely to an
${\cal O}_{K}$-valued point
$\underline{y}$
of
${\bf P}^{N}_{{\cal O}_{K}}$.
Let
$\overline{{\cal O} (1)}$
be the hermitian line bundle on
${\bf P}^{N}_{{\cal O}_{K}}$,
where the hermitian metrics at the infinite places are given by
$$\left\| x_{0} \right\| := \left( 1 + \left| \frac{x_{1}}{x_{0}} \right|^{2} +
\ldots + \left| \frac{x_{N}}{x_{0}} \right|^{2} \right)^{- \frac{1}{2}}
{}~~~~~~~~~ ({\rm i.e.~~} \left\| x_{i} \right\| := \left( \left|
\frac{x_{0}}{x_{i}} \right|^{2} + \ldots + 1 + \ldots + \left|
\frac{x_{N}}{x_{i}} \right|^{2} \right)^{- \frac{1}{2}} ) ~~.$$
Then
$h = h_{\overline{{\cal O} (1)}}$
is the height defined by
$\overline{{\cal O} (1)}$
in the sense of [BoGS, Definition 3.1.; \linebreak formula (3.1.6)].
}
\end{thm}

\begin{thm}
{\bf Remark.}
{\rm
We need a better understanding of
${\cal O} (j^{*} (\Theta))$.
By Riemann's Theorem [GH, Chapter 2, \S 7] one has
$\Theta = \frac{1}{(g-1)!} j_{*} ((x) \times C^{g-1})$,
where
$j: C^{g} \stackrel{p}{\longrightarrow} C^{(g)} \stackrel{c}{\longrightarrow}
J$
factors into a morphism finite flat of degree
$g!$
and a birational morphism. So
$j^{*} (\Theta)$
is an effective divisor containing the summands
$\pi_{k}^{*} (x)$,
where
$\pi_{k}: C^{g} \longrightarrow C$
denotes
$k$-th
projection.
$${\cal O} \Big( j^{*} \left( \Theta \right) \Big) = \bigotimes_{k=1}^{g}
\pi_{k}^{*} \Big( {\cal O} (x) \Big) \otimes {\cal O} \Big( p^{*} (R) \Big)$$
Intuitively, the divisor
$R$
on
$C^{(g)}$
corresponds to the divisors on
$C$
moving in a linear system. This can be made precise, but we will not need that
here.
}
\end{thm}

\begin{thm}
{\bf Remark.}
{\rm It is a difficulty that there are no regular projective models available
for
$J$
and
$C^{g}$,
such that arithmetic intersection theory does not work immediately. So we
follow [BoGS, Remark after Proposition 3.2.1.] and consider a projective (not
necessarily regular) model of
$C^{g}$,
namely
${\cal C}^{g} := {\cal C} \times_{{\rm Spec} {\cal O}_{K}} \ldots \times_{{\rm
Spec} {\cal O}_{K}} {\cal C}$.
Hereon let
$\overline{\cal T}$
be a line bundle extending
$\bigotimes_{k=1}^{g} \pi_{k}^{*} {\cal O} (x)$
equipped with a hermitian metric. One has to define a height
$h_{\overline{\cal T}}$
induced by
$\overline{\cal T}$.

Consider more generally a projective (singular) arithmetic variety
${\cal X} / {\cal O}_{K}$
and a hermitian line bundle
$\overline{\cal U}$
on
$\cal X$.
Then there is a morphism
$\iota: {\cal X} \longrightarrow P$
into a projective variety
$P$
smooth over
${\rm Spec} ~ {\cal O}_{K}$
and a line bundle
${\cal U_{P}}$
on
$P$
such that
$\iota^{*} (\cal U_{P}) = {\cal U}$
(see [Fu, Lemma 3.2.], cf. [BoGS, Remark 2.3.1.ii)]). We can even choose
$\iota$
in such a way that the hermitian metric on
$\overline{\cal{U}}$
is a pullback of one on
${\cal U}_{P}$
(e.g. as a closed embedding).
$$\iota^{*} \Big( \overline{{\cal U}_{P}} \Big) = \overline{{\cal U}}$$
Then for an
${\cal O}_{K}$-valued
point
$\underline{y}$
of
${\cal X}$
one defines
\begin{eqnarray*}
h_{\overline{\cal U}} \Big( \underline{y} \Big) & := & h_{\overline{\cal
U}_{P}} \Big( \iota_{*} (\underline{y}) \Big) \\
 & = & \stackrel{\wedge}{\deg} \Big( \stackrel{\wedge}{c_{1}} (\overline{{\cal
U}_{P}} ) \Big| \iota_{*} (y) \Big) ~~,
\end{eqnarray*}
where
$( . | . )$
denotes the pairing
$\stackrel{\wedge}{{\rm CH}^{1}} (P) \times {\rm Z}_{1} (P) \longrightarrow ~
\stackrel{\wedge}{{\rm CH}^{1}} ({\rm Spec} ~ {\cal O}_{K} )_{\Bbb{Q}}$
from [BoGS, 2.3.]. In [BoGS, Remark after Proposition 3.2.1.] independence of
this definition of the
$\iota$
chosen is shown. In particular it becomes clear at this point that the pairing
$\Big( \stackrel{\wedge}{c_{1}} ( . ) \Big| . \Big)$
can be extended to arbitrary (singular) projective arithmetic varieties over
${\cal O}_{K}$
and satisfies the projection formula
$$ \Big( \stackrel{\wedge}{c_{1}} (L) \Big| f_{*} (Z) \Big) = \Big(
\stackrel{\wedge}{c_{1}} (f^{*} (L)) \Big| Z \Big) ~~ .$$
}
\end{thm}

\begin{thm}
\label{sing}
{\bf Remark.}
{\rm If
${\cal X} / {\cal O}_{K}$
is a regular arithmetic variety, one has another pairing
\begin{eqnarray*}
[.,.]: \stackrel{\wedge}{{\rm CH}^{1}} ({\cal X}) ~ \times
\stackrel{~^{\scriptstyle{\wedge}}}{{\rm CH}_{1}} ({\cal X}) & \longrightarrow
& \stackrel{\wedge}{{\rm CH}^{1}} ({\rm Spec} ~ {\cal O}_{K})_{\Bbb Q} \\
 (z, y) & \mapsto & \pi_{*} [z \cdot y] ~~.
\end{eqnarray*}
We note that also
$\Big[ \stackrel{\wedge}{c_{1}} (.) , . \Big]$
can be extended to arbitrary (singular away from the generic fiber) projective
arithmetic varieties. One has to represent
$y$
by a cycle
$(Y, g_{Y})$
and to put
$$\Big[ \stackrel{\wedge}{c_{1}} (\overline{\cal U}) , Y \Big] := \Big(
\stackrel{\wedge}{c_{1}} (\overline{\cal U}) \Big| Y \Big) + \Big(0, \Big(
\int_{{\cal X} ({\Bbb C})} ~ g_{Y} \omega_{\stackrel{\wedge}{c_{1}}
(\overline{\cal U})} \Big)_{\sigma: K \hookrightarrow {\Bbb C}} \Big)$$
obtaining a pairing satisfying the projection formula
$$\Big[ \stackrel{\wedge}{c_{1}} (f^{*} (\overline{\cal U})) , y \Big] = \Big[
\stackrel{\wedge}{c_{1}} (\overline{\cal U}) , f_{*} (y) \Big]$$
for
$f$
proper and smooth on the generic fiber. In particular, independence of the
cycle chosen carries over from the regular case. Indeed, concerning a trivial
arithmetic
$1$-cycle
one is automatically reduced to surfaces and resolution of singularities is
known for two-dimensional schemes [CS, Chapter XI by M. Artin]. Let
$f$
be one. Note that for cycles with
$\omega_{y} = 0$
the push-forward
$f_{*}$
makes sense for any proper
$f$.

If
$f$
is a proper birational map inducing an isomorphism on the generic fibers one
has
$f_{*} f^{*} w = w$
for arithmetic one-cycles and therefore
\begin{equation}
\label{fun}
\Big[ \stackrel{\wedge}{c_{1}} (f^{*} (\overline{\cal U})) , f^{*} (w) \Big] =
\Big[ \stackrel{\wedge}{c_{1}} (\overline{\cal U}) , w \Big] ~~.
\end{equation}
This is useful for the special case of a (singular) projective arithmetic
surface. There
$[ . , . ]$
can be specialized to a pairing
$\Big[ \stackrel{\wedge}{c_{1}} (.) , \stackrel{\wedge}{c_{1}} (.) \Big]$
between hermitian line bundles. This one is symmetric. Indeed, formula
(\ref{fun}) tells us, that it is enough to show that after pullback by a
birational morphism. But for regular arithmetic surfaces symmetry is clear.
}
\end{thm}

\begin{thm}
\label{class}
{\bf Lemma.}
Let
${\cal X} / {\cal O}_{K}$
be a (singular) projective arithmetic variety and
$X/K$
its generic fiber which is assumed to be regular. Further, let
$D$
be a divisor on
$X$
and
$\overline{\cal U}$
be a hermitian line bundle extending
${\cal O} (D)$.
Then
$h_{D} = h_{\overline{\cal U}} + {\rm O} (1)$
for
$K$-valued
points of
$X.$
\newline
{\bf Proof.}
{\rm There is a very ample line bundle on
$X$
that can be extended to
${\cal X}$.
So we may assume
$D$
to be basepoint-free (very ample). Then the two height functions arise from the
situations
$$\begin{array}{cccccccccccccc}
      &    &    &    & \overline{{\cal O} (1)} & ~~ & & ~~ & & &
{}~~~~\overline{\cal U}~~~~ & & \overline{{\cal U}_{P}}~~ & \\
   &     &    &     & |  & ~~ & {\rm and} & ~~ & & & | & & | & \\
{\rm Spec} ~ {\cal O}_{K} & \stackrel{y}{\hookrightarrow} & {\cal X} & -
\stackrel{i}{-} \rightarrow & {\bf P}^{N}_{{\cal O}_{K}} & ~~ & & ~~ & {\rm
Spec} ~ {\cal O}_{K} & \stackrel{y}{\hookrightarrow} & {\cal X} &
\stackrel{\iota}{\longrightarrow} & P & ~~.
\end{array}$$
Here
$i$
is the rational map defined by an extension
${\cal U}^{'}$
of
${\cal O} (D)$
over
${\cal X}$.
In the generic fiber
$i$
is defined everywhere. Note that
$iy$
is a morphism by the valuative criterion. Of course, it comes from sections of
the line bundle
$y^{*} {\cal U}^{\prime}$.
Note that
${\cal U}^{'}$
is equipped with a hermitian metric induced by that on
${\cal O} (1)$.
$\iota: {\cal X} \longrightarrow P$
is a morphism into a smooth scheme as described above. Thus
\begin{eqnarray*}
h_{D} (y) & = & \stackrel{\wedge}{\deg} \left( \stackrel{\wedge}{c_{1}} \left(
\overline{{\cal O} (1) } \right) \Big| (iy)_{*} ({\rm Spec} ~ {\cal O}_{K})
\right) \\
 & = & \stackrel{\wedge}{\deg} \left( \stackrel{\wedge}{c_{1}} \left( (iy)^{*}
\overline{{\cal O} (1) } \right) \Big| {\rm Spec} ~ {\cal O}_{K} \right) \\
 & = & \stackrel{\wedge}{\deg} \left( \stackrel{\wedge}{c_{1}} \left( (iy)^{*}
\overline{{\cal O} (1) } \right) \right) \\
 & = & \stackrel{\wedge}{\deg} \left( \stackrel{\wedge}{c_{1}} \left( y^{*}
\left( \overline{{\cal U}^{'}} \right) \right) \right)
\end{eqnarray*}
and, correspondingly,
\begin{eqnarray*}
h_{\overline{\cal T}} (y) & = & \stackrel{\wedge}{\deg} \left(
\stackrel{\wedge}{c_{1}} \left( \overline{{\cal U}_{P}} \right) \Big| (\iota
y)_{*} ({\rm Spec} ~ {\cal O}_{K}) \right) \\
 & = & \stackrel{\wedge}{\deg} \left( \stackrel{\wedge}{c_{1}} \left( y^{*}
\left( \overline{\cal U } \right) \right) \Big| {\rm Spec} ~ {\cal O}_{K}
\right) \\
 & = & \stackrel{\wedge}{\deg} \left( \stackrel{\wedge}{c_{1}} \left( y^{*}
\left( \overline{\cal U } \right) \right) \right) ~~~ .
\end{eqnarray*}
But
$\overline{{\cal U}^{'}}$
and
$\overline{\cal U}$
coincide as line bundles on the generic fiber. As bundles their difference is
some
${\cal O} (E)$
where
$E$
is a divisor contained in the special fibers of
${\cal X}$,
while the hermitian metrics differ by a continuous, hence bounded, factor.
Therefore, the first arithmetic Chern classes of the pullbacks considered
differ only at the infinite and a finite number of finite places by bounded
summands.
}
\begin{center}
$\Box$
\end{center}
\end{thm}

\begin{thm}
{\bf Remark.} {\rm a) When one considers
$L$-valued
points instead of
$K$-valued
ones, where
$L$
is a number field with
$[L:K] = d$,
the error term becomes
${\rm O}(d)$;
i. e. there is a constant
$C$
such that
$$ \Big| h_{D} (x) - h_{\overline{\cal U}} (x) \Big| < C \cdot d$$
for
$L$-valued
points
$x$
of
$X$
and an arbitrary number field
$L/K$.
The reason for that is simply that the number of the critical places occuring
grows as
${\rm O}(d)$.
\newline
b) The lemma can be applied to
${\cal X} = {\cal C}^{g}$
and
$D = \sum_{k=1}^{g} \pi_{k}^{*} (x)$,
since
$\bigotimes_{k=1}^{g} \pi_{k}^{*} {\cal O} (\overline{x})$
extends
${\cal O} (D)$.
}
\end{thm}

\begin{thm}
\label{div}
{\rm The height defined by an extension
${\cal U}$
of
$\bigotimes_{k=1}^{g} \pi_{k}^{*} ({\cal O} (x))$
is understood by the following
}
\newline
{\bf Proposition.}
On
${\cal C} / {\cal O}_{K} $
let
$\overline{\cal S}$
be the line bundle
${\cal O} (\overline{x})$,
where
$\overline{x}$
denotes the closure of
$x$
in
$\cal C$,
equipped with a hermitian metric. For
$L$-valued
points
$P = (P_{1}, \ldots ,P_{g})$
of
$C^{g}$
we consider the divisor
$\underline{P} := (P_{1}) + \ldots + (P_{g})$
on
$C$.
Then
$$ h_{\overline{\cal S}} (\underline{P}) = h_{\overline{\cal U}} (P) + O(d)
{}~~.$$
\end{thm}
{\bf Proof.}
By [BoGS, Proposition 3.2.2.ii)] we may assume
${\cal U} = {\cal O} \Big( \sum_{k=1}^{g} \pi_{k}^{*} (\overline{x}) \Big)$,
where
$\overline{x}$
is the closure of
$x$
in
${\cal C}$.
The extensions of
$P$
and
$\underline{P}$
over
$\cal C$
and
${\cal C}^{g}$
will be denoted by
$(\overline{P_{1}}) + \ldots + (\overline{P_{g}})$,
respectively
$(\overline{P_{1}}, \ldots , \overline{P_{g}})$.
Then
\begin{eqnarray*}
h_{\overline{\cal S}} \Big( (\overline{P_{1}}) + \ldots + (\overline{P_{g}})
\Big) & = & \stackrel{\wedge}{\deg} \Big( \stackrel{\wedge}{c_{1}}
(\overline{\cal S}) \Big| (\overline{P_{1}}) + \ldots + (\overline{P_{g}})
\Big) \\
 & = & \sum_{k=1}^{g} \stackrel{\wedge}{\deg} \Big( \stackrel{\wedge}{c_{1}}
(\overline{\cal S}) \Big| (\overline{P_{k}} ) \Big) \\
 & = & \sum_{k=1}^{g} \stackrel{\wedge}{\deg} \Big( \stackrel{\wedge}{c_{1}}
\Big( \pi_{k}^{*} (\overline{\cal S}) \Big) \Big| (\overline{P_{1}}, ~ \ldots ~
, \overline{P_{g}} ) \Big) ~~~~~~~~~~~{\rm "projection ~formula"}\\
 & = & \stackrel{\wedge}{\deg} \left( \stackrel{\wedge}{c_{1}} \left(
\bigotimes_{k=1}^{g} \pi_{k}^{*} \left( \overline{\cal S} \right) \right)
\bigg| \left( \overline{P_{1}}, ~ \ldots ~ , \overline{P_{g}} \right) \right)
{}~~.
\end{eqnarray*}
But by construction
$\bigotimes_{k=1}^{g} \pi_{k}^{*} (\overline{\cal S})$
is the line bundle
${\cal U}$,
equipped with a hermitian metric (and by definition the formula
$\Big( \stackrel{\wedge}{c_{1}} (\bigotimes_{k} \overline{{\cal L}_{k}}) \Big|
Z \Big) = {\displaystyle\sum}_{k} \Big( \stackrel{\wedge}{c_{1}}
(\overline{{\cal L}_{k}}) \Big| Z \Big) $
holds in singular case, too). So we have
$$h_{\overline{\cal S}} \Big( (\overline{P_{1}}) + \ldots + (\overline{P_{g}})
\Big) = h_{\overline{{\cal U}}^{'}} \Big( (\overline{P_{1}} , ~ \ldots ~ ,
\overline{P_{g}}) \Big) ~~~,$$
\newpage
{}~
\newline
where
$\overline{{\cal U}}^{'}$
differs from
$\overline{\cal U}$
only by the hermitian metric. The claim follows from [BoGS, Proposition
3.2.2.i)].
\begin{center}
$\Box$
\end{center}

\begin{thm}
\label{two}
{\bf Corollary.}
Let
$P \in C^{g}(L)$
and
$\underline{P}$
be the associated divisor on
$C$.
Then
$$h_{\Theta} ({\cal O} (\underline{P})) = h_{\overline{\cal S}} (\underline{P})
+ h_{R} (p_{*} P) + {\rm O}(d) ~~,$$
where
$h_{R}$
denotes the height for points of
$C^{(g)}$
defined by
$R$.
\end{thm}

\section{An observation concerning the tautological line \newline bundle}
In this section we start analyzing the fundamental definition 1.3. First we
will consider only varieties over number fields and forget about integral
models.

\begin{thm}
{\bf Definition.}
Let
$\Delta$
be the diagonal in
$C \times C$.
Then
$$\underline{\underline{\cal E}} := \bigotimes_{k=1}^{g} \pi_{k,g+1}^{*} \Big(
{\cal O} (\Delta) \Big) $$
will be called the {\rm tautological line bundle} on
$C^{g} \times C$.
Note that the restriction of
$\underline{\underline{\cal E}}$
to
$\{ (P_{1}, \ldots , P_{g}) \} \times C$
equals
${\cal O} (P_{1} + \ldots + P_{g})$.
By construction
$\underline{\underline{\cal E}}$
is the pullback of some line bundle
$\cal E$,
said to be {\rm the tautological} one on
$C^{(g)} \times C$.
$$\underline{\underline{\cal E}} = (p \times id)^{*} ({\cal E})$$
\end{thm}

\begin{thm}
{\bf Proposition.}
We have
$\det R \pi_{*} {\cal E} = {\cal O}_{C^{(g)}} (-R)$.
\end{thm}

\begin{thm}
{\rm This will be a direct consequence of the following} \newline
{\bf Lemma.}
Let
$\underline{\cal E} := {\cal E} \otimes \pi_{C}^{*} ({\cal O} (-x))$ be a
tautological line bundle fiber-by-fiber of degree
$g-1$.
Then
$$\det R \pi_{*} \underline{\cal E} = {\cal O}_{C^{(g)}} \Big( -c^{*} (\Theta)
\Big) ~~.$$
\end{thm}
{\bf Proof.}
The canonical map
$c: C^{(g)} \longrightarrow J$
is given by
${\cal E}$
using Picard functoriality. So for a tautological line bundle
${\cal M}$,
fiber-by-fiber of degree
$g$
on
$J \times C$,
one has
$${\cal E} = (c \times id)^{*} {\cal M} \otimes \pi^{*} {\cal H} ~~,$$
where
$\cal H$
is a line bundle on
$C^{(g)}$.
Putting
${\cal M}_{0} := {\cal M} \otimes \pi_{C}^{*} {\cal O} (-x)$,
where
$\pi_{C}: J \times C \longrightarrow C$
denotes here the canonical projection from
$J \times C$,
we get
$$\underline{\cal E} = (c \times id)^{*} {\cal M}_{0} \otimes \pi^{*} {\cal H}
{}~~.$$
It follows
\begin{eqnarray*}
\det R \pi_{*} \underline{\cal E} & \cong & \det R \pi_{*} \Big[ (c \times
id)^{*} {\cal M}_{0} \otimes \pi^{*} {\cal H} \Big] \\  & = & \det R \pi_{*}
\Big[ (c \times id)^{*} {\cal M}_{0} \Big] \\ & = & c^{*} \det R \pi_{*} {\cal
M}_{0} ~~,
\end{eqnarray*}
where we first used the projection formula, which is particularly simple here,
since line bundles, fiber-by-fiber of degree
$g-1$,
have relative Euler characteristic
$0$,
and afterwords noted \pagebreak
that the determinant of cohomology commutes with arbitrary base change
\cite{Knudsen/Mumford 76}. But by [MB, Proposition 2.4.2] or [Fa, p. 396] we
know
$\det R \pi_{*} {\cal M}_{0} = {\cal O}_{J} (- \Theta)$.
The assertion follows.
\begin{center}
$\Box$
\end{center}

\begin{thm}
{\bf Proof of the Proposition.}
{\rm The short exact sequence
$$0 \longrightarrow \underline{\cal E} \longrightarrow {\cal E} \longrightarrow
{\cal E}|_{C^{(g)} \times \{ x \} } \longrightarrow 0 $$
gives
\begin{eqnarray*}
\det R \pi_{*} {\cal E} & = & \det R \pi_{*} \underline{\cal E} \otimes {\cal
O} \left( {\frac{\scriptstyle 1}{\scriptstyle g!}} p_{*} \left( \sum_{k=1}^{g}
\pi_{k}^{*} (x) \right) \right) \\
  & = & {\cal O} \left( - c^{*} (\Theta) \right) \otimes {\cal O} \left(
\frac{\scriptstyle 1}{\scriptstyle g!} p_{*} \left( \sum_{k=1}^{g} \pi_{k}^{*}
(x) \right) \right) \\
  & = & {\cal O} (-R) ~~.
\end{eqnarray*}
\begin{center}
$\Box$
\end{center}
}
\end{thm}

\begin{thm}
{\bf Corollary.}
\label{glatt}
$\det R \pi_{*} ({\cal E} \otimes \pi^{*} {\cal O} (R)) = {\cal O}_{C^{(g)}}$.
\newline
{\bf Proof.}
{\rm This is the projection formula for the determinant of cohomology.}
\begin{center}
$\Box$
\end{center}
\end{thm}

\begin{thm}
{\rm Let
$\cal J$
be the N\'eron model of the Jacobian
$J$
of
$C$.
It is smooth over
${\cal O}_{K}$,
consequently
${\cal J} \times_{{\rm Spec} {\cal O}_{K}} {\cal C}$
is smooth over
${\cal C}$
and therefore regular. We note that any
$K$-valued
point of
$J$
can be extended uniquely to an
${\cal O}_{K}$-valued
point of
$\cal J$.

On
$J \times C$
we have a tautological line bundle
${\cal M}$,
fiber-by-fiber of degree
$g$.
${\cal M}$ can be extended over
${\cal J} \times_{{\rm Spec} {\cal O}_{K}} {\cal C}$.
For this let
${\cal M} = {\cal O} (D)$
with some Weil divisor
$D$
on
$J \times C$.
Its closure
$\overline{D}$
in
${\cal J} \times_{{\rm Spec} {\cal O}_{K}} {\cal C}$
is obviously flat over
${\cal O}_{K}$
and therefore it has codimension
$1$.
We choose the extension
${\cal O} ({\overline{D}})$
and denote it by
$\cal M$
again.

$\cal M$
is a perfect complex of
${\cal O}_{{\cal J} \times_{{\rm Spec} {\cal O}_{K}} {\cal C}}$-modules. For
the existence of the Knudsen-Mumford determinant we need that
$$\pi: {\cal J} \times_{{\rm Spec} {\cal O}_{K}} {\cal C} \longrightarrow {\cal
J}$$
has finite
${\rm Tor}$-dimension.
For this there exists a closed embedding
${\cal C} \longrightarrow P$,
where
$P$
is smooth over
${\cal O}_{K}$.
Thus
$\pi$
factorizes as
$${\cal J} \times_{{\rm Spec} {\cal O}_{K}} {\cal C}
\stackrel{i}{\hookrightarrow} {\cal J} \times_{{\rm Spec} {\cal O}_{K}} P
\stackrel{{\rm smooth}}{\twoheadrightarrow} {\cal J} ~~.$$
By [SGA 6, Expos\'e III, Proposition 3.6] it is enough to show that
$i$
has finite
${\rm Tor}$-dimension.
But
$i_{*}$
is exact and
${\cal J} \times_{{\rm Spec} {\cal O}_{K}} P$
is regular implying quasi-coherent sheaves have locally finite free resolutions
of bounded length.

${\cal M}$
has, relative to
$\pi$,
Euler characteristic
$1$.
Therefore
$\cal M$
can be changed by an inverse image of a line bundle on
${\cal J}$,
trivial on the generic fiber, in such a way, that we are allowed to assume
$$\det R \pi_{*} {\cal M} \cong {\cal O}_{\cal J} ~~.$$
}
\end{thm}

\section{Choosing hermitian metrics continuously depending on moduli space}

\setcounter{equation}{1}
\begin{thm}
{\bf Fact.}
On
${\cal M}_{\Bbb{C}}$
there exists a hermitian metric
$\underline{h}$
such that for every point
$y \in J ({\Bbb{C}})$
the curvature form satisfies
$$c_{1} ({\cal M}_{{\Bbb{C}},y}, \underline{h}_{y}) = g \omega$$
on
$( \{ y \} \times C) ({\Bbb{C}}) \cong C (\Bbb{C})$.
\newline
{\bf Proof.}
{\rm The statement is local in
$C^{g} (\Bbb{C})$
by partition of unity. By the Theorem on cohomology and base change
$R^{0} \pi_{*} {\cal M}_{\Bbb{C}} (g-1)$
is locally free and commutes with arbitrary base change. Hence there exists,
locally on
$J (\Bbb{C})$,
a rational section
$s$
of
$\cal M$
that is neither undefined nor identically zero in any fiber.

First we choose an arbitrary hermitian metric
$\| . \|$
on
${\cal M}_{\Bbb{C}}$.
Then
\begin{equation}
\omega^{'} := -d_{C} d_{C}^{c} \log \| s \|^{2}
\end{equation}
defines a smooth
$(1,1)$-form
on
$(J \times C) ({\Bbb{C}}) \backslash {\rm div} (s)$,
that is fiber-by-fiber the curvature form to be considered. Since construction
(2) is independent of
$s$
as soon as it makes sense at a point, we obtain
$\omega^{'}$
as a smooth
$(1,1)$-form
on
$(J \times C) (\Bbb{C})$
closed under
$d_{C}$
and cohomologous to
$g \omega$
on
$\{ y \} \times C(\Bbb{C})$
for any
$y \in C^{g} (\Bbb{C})$.

The setup
$\| . \|_{\underline{h}} = f \cdot \| . \|$
gives the equation
\begin{equation}
\omega^{'} - g \omega = d_{C} d_{C}^{c} \log | f |^{2} ~~.
\end{equation}
But
$d d^{c}$
is an elliptic differential operator on the Riemann surface
$C (\Bbb{C})$,
so by Hodge theory it permits a Green`s operator
$G$
compact with respect to every Sobolew norm
$\| . \|_{\alpha}$.
Consequently, there exists a solution
$f$
of (3) being smooth on
$(J \times C) (\Bbb{C})$.
}
\begin{center}
$\Box$
\end{center}
\end{thm}

\begin{thm}
{\rm We note, that
$\det R \pi_{*} {\cal M} \cong {\cal O}_{\cal J}$
and the isomorphism is uniquely determined up to units of
${\cal O}_{K}$.
Namely, one has
${\rm Aut}_{{\cal O}_{\cal J}} ({\cal O}_{\cal J}) = \Gamma ({\cal J}, {\cal
O}_{\cal J}^{*})$
and already
$\Gamma (J, {\cal O}_{J}^{*})$
consists of constants only. In particular, there is a unitary section, uniquely
determined up to units of
${\cal O}_{K}$,
$${\bf 1} \in \Gamma ({\cal J}, \det R \pi_{*} {\cal M}) ~~.$$
}
\end{thm}

\begin{thm}
{\bf Corollary.}
Let
$R \in \Bbb{R}$.
Then, on
${\cal M}_{\Bbb{C}}$
there exists exactly one hermitian metric
$h$,
such that for every point
$y \in J(\Bbb{C})$
the curvature form
$c_{1} ({\cal M}_{{\Bbb{C}}, y}, h_{y}) = \omega$
and for the Quillen metric one has
$$h_{Q,h} ({\bf 1}) = R ~~,$$
where
${\bf 1} \in \Gamma (\{ y \} , \det R \pi_{*} {\cal M}_{{\Bbb{C}}, y} )$.
\newline
{\bf Proof.}
{\rm Let
$\underline{h}$
be the hermitian metric from the preceeding fact. We may replace
$\underline{h}$
by
$f \cdot \underline{h}$
with
$f \in C^{\infty} (J(\Bbb{C}))$
without any effort on the curvature forms, since they are invariant under
scalation. As
$\cal M$
has relative Euler characteristic
$1$,
exactly
$$h := \frac{R}{h_{Q, \underline{h}} ({\bf 1})} \cdot \underline{h} $$
satisfies the conditions required.}
\begin{center}
$\Box$
\end{center}
\end{thm}

\begin{thm}
{\rm We have to consider
${\cal M}_{\Bbb{C}}$
on
$J({\Bbb{C}}) \times (\coprod_{\sigma: K \hookrightarrow {\Bbb{C}}} C
(\Bbb{C}))$.
The metric
$h$
on
${\cal M}_{\Bbb{C}}$
has to be invariant under
$F_{\infty}$,
its curvature form is required to be
$g \omega$
and we want to realise
\begin{equation}
\prod_{\sigma: K \hookrightarrow \Bbb{C}} h_{Q,h} ({\bf 1}) = 1
\end{equation}
simultaneously for all
$y \in J(\Bbb{C})$.

The first is possible since
$\omega$
is invariant under
$F_{\infty}$
and the corollary above already gives conditions uniquely determining
$h$.
(4) can be obtained by scalation with a constant factor over all
$J({\Bbb{C}}) \times (\coprod_{\sigma: K \hookrightarrow {\Bbb{C}}}
C({\Bbb{C}}))$.

Altogether, for every line bundle of degree
$g$
on
$C$
we have found a distinguished hermitian metric and seen that it depends, in
some sense, continuously on the moduli space $J$.
One obtains
}
\end{thm}

\begin{thm}
{\bf Proposition.}
\label{Ara}
Let
$K$
be a number field and
$({\cal C} / {\cal O}_{K}, \omega)$
a regular connected Arakelov surface. Then, on the (non proper) Arakelov
variety
$({\cal J} \times_{{\rm Spec} {\cal O}_{K}} {\cal C}, \pi_{C}^{*} \omega)$
there is a hermitian line bundle
$\overline{\cal M}$
with the following properties.
\newline
{\rm a)}
$$(c \times id)^{*} ({\cal M} |_{J \times C}) = {\cal E} \otimes \pi^{*} {\cal
O} (R)$$
is the modified tautological line bundle found in section 3.
\newline
{\rm b)} The hermitian metric
$h$
on
${\cal M} |_{\coprod_{\sigma: K \hookrightarrow \Bbb{C}} (J \times C)
({\Bbb{C})}}$
is invariant under
$F_{\infty}$
and has curvature form
$g \omega$.
\newline
{\rm c)} For any
$y \in J(K)$
one has
$\overline{\{ y \}} \subseteq {\cal J}$
and
$$\stackrel{\wedge}{\deg} \left( \det R \pi_{*} \left( {\cal M} |_{\overline{\{
y \} } \times_{{\rm Spec} {\cal O}_{K}} {\cal C}}, h_{{\cal M},y} \right) , \|
. \|_{Q,h} \right) = 0 ~~.$$
{\bf Proof.}
{\rm b) and c) are clear. For a) we know
${\cal E} = (c \times id)^{*} ({\cal M}|_{J \times C}) \otimes \pi^{*} {\cal
H}$
from 4.3. But
$\cal H$
is determined by
$\det R \pi_{*} {\cal E} = {\cal O}_{C^{(g)}} (-R)$
and
$\det R \pi_{*} ({\cal M}|_{J \times C}) = {\cal O}_{J}$.
\begin{center}
$\Box$
\end{center}
}
\end{thm}

\section{An integral model of the symmetric power}
\begin{thm}
{\bf Remark.}
{\rm It turns out here to be very inconvenient to work directly with the
N\'eron model
$\cal J$
of the Jacobian of
$C$.
When one considers the tautological line bundle
${\cal M} |_{J \times C}$,
fiber-by-fiber of degree
$g$
on
$J \times C$
with
$\det R \pi_{*} {\cal M}|_{J \times C} \cong {\cal O}_{J},$
then
${\cal M}|_{J \times C}$
will even have an (up to constant factor) canonical section.
$$\pi_{*} {\cal M}|_{J \times C} \cong {\cal O}_{J}$$
But this section is zero over a codimension two subset of
$J$
such that one is led to blow up this subset.
}
\end{thm}

\begin{thm}
{\bf Lemma.}
{\rm a)}
$C^{(g)}$
is a projective variety.
\newline
{\rm b)}
The divisor
$S := \frac{1}{g!} p_{*} \Big( \sum_{k=1}^{g} \pi_{k}^{*} (x) \Big) =
\frac{1}{(g-1)!} p_{*} \Big( (x) \times C^{g-1} \Big)$
"one of the points is
$x$"
on
$C^{(g)}$
is ample.
\newline
{\bf Proof.}
{\rm
a) There are at least two good reasons for that. First
$C^{(g)}$
is proper as a quotient of the proper variety
$C^{g}$
and b) gives an ample line bundle. On the other hand we can give a high-powered
argument as follows.
$C^{(g)}$
is the Hilbert scheme
${\rm Hilb}_{C/K}^{g}$
by [CS, Chapter VII by J. S. Milne, Theorem 3.13] and this is known to be
projective for a long time [FGA, Expos\'e 221, Theorem 3.2].
\pagebreak
\newline
b) By [EGA III, Proposition 2.6.2] it is enough to show that
$p^{*} {\cal O} (S) = \bigotimes_{k=1}^{g} \pi_{k}^{*} {\cal O} (x)$
is an ample line bundle on
$C^{g}$,
which is obvious.
}
\begin{center}
$\Box$
\end{center}
\end{thm}

\begin{thm}
{\bf Proposition.}
The morphism
$c: C^{(g)} \longrightarrow J$
is the blow-up of some ideal sheaf
${\cal I} \subseteq {\cal O}_{J}$
with
$$c^{-1} {\cal I} = {\cal O} (-NR) ~~,$$
where
$N$
is a positive integer.
\newline
{\bf Proof.}
{\rm
$c$
is birational and by the lemma it is a projective morphism. So it is a blow-up
of some ideal sheaf
${\cal I} \subseteq {\cal O}_{J}$.
Going through the lines of the proof of [Ha, Chapter II, Theorem 7.17] one sees
that
$c_{*} ({\cal O} (NS))$
for
$N \gg 0$,
up to tensor product with line bundles in order to make them ideal sheaves, can
be used as
$\cal I$.
But
$c^{*} {\cal O} (\Theta) = {\cal O} (S) \otimes {\cal O} (R)$
gives
\begin{eqnarray*}
c_{*} ({\cal O} (NS)) & = & c_{*} \Big( {\cal O} (-NR) \otimes c^{*} {\cal O}
(N \Theta) \Big) \\
  & = & c_{*} \Big( {\cal O} (-NR) \Big) \otimes {\cal O} (N \Theta)
\end{eqnarray*}
and therefore
${\cal I} = c_{*} {\cal O} (-NR)$
for some
$N \gg 0$.

We have a short exact sequence
$$0 \longrightarrow {\cal O}_{C^{(g)}} (-NR) \longrightarrow {\cal O}_{C^{(g)}}
\longrightarrow {\cal O}_{R_{N}} \longrightarrow 0 ~~,$$
where
$R_{N}$
denotes the
$N$-th
infinitesimal neighbourhood of
$R$.
It follows exactness of
$$0 \longrightarrow c_{*} {\cal O}_{C^{(g)}} (-NR) \longrightarrow {\cal O}_{J}
\longrightarrow c_{*} {\cal O}_{R_{N}} ~~.$$
Now the image of
${\cal O}_{J} \longrightarrow c_{*} {\cal O}_{R_{N}}$
is the structure sheaf of the scheme-theoretic image
$I_{N}$
of
$R_{N}$
in
$J$.
So
${\cal I} = c_{*} {\cal O}_{C^{(g)}} (-NR) = {\cal I}_{I_{N}} \subseteq {\cal
O}_{J}$.
But
${\cal O}_{C^{(g)}} / c^{-1} {\cal I}_{I_{N}} = c^{*} ({\cal O}_{C^{g}} / {\cal
I}_{I_{N}})$
and therefore
$c^{-1} {\cal I} = c^{-1} {\cal I}_{I_{N}}$
is the ideal sheaf of
$I_{N} \times_{J} C^{(g)}$
in
${\cal O}_{C^{(g)}}$.

$c^{-1} {\cal I}$
is known to be invertible, so
$I_{N} \times_{J} C^{(g)}$
is necessarily pure of codimension
$1$
and by construction it contains the scheme
$R_{N}$.
But
$R = c^{*} c_{*} S - S$
contains with one point its complete fiber in
$c: C^{(g)} \longrightarrow J$.
So
$I_{N} \times_{J} C^{(g)}$
must be an infinitesimal thickening of
$R_{N}$.
On the other hand, when one replaces
$R$
by
$S$
and considers the scheme-theoretic image
$\underline{I_{N}}$
of the
$N$-th
infinitesimal neighbourhood
$S_{N}$,
then
$\underline{I_{N}} \times_{J} C^{(g)} \supseteq I_{N} \times_{J} C^{(g)}$
is a pure codimension
$1$
subscheme not containing any thickening of
$R_{N}$,
but only other additional summands (it corresponds to the divisor
$c^{*} (N \Theta)$).
So, necessarily
$I_{N} \times_{J} C^{(g)} = R_{N}$
and
$$c^{-1} {\cal I} = {\cal O} (-NR) ~~.$$
}
\begin{center}
$\Box$
\end{center}
\end{thm}

\begin{thm}
{\rm Denote by
$\tilde{\cal J}$
the normalization of the blow-up of
$\cal J$
with respect to some extension
$\underline{\cal I}$
of the ideal sheaf
$\cal I$
over
$\cal J$.
}
\newline
{\bf Facts.}
{\rm a)}
$\tilde{\cal J}$
is some (singular) arithmetic variety proper over
$\cal J$.
\newline
{\rm b)}
It is an integral model of
$C^{(g)}$.
\newline
{\rm c)}
On
$\tilde{\cal J}$
one has the line bundle
$${\cal R} := (c^{-1} \underline{\cal I})^{\vee}$$
extending
${\cal O} (NR)$
for some
$N > 0$.
\newline
{\rm
Note that we do not know whether we have an extension of
${\cal O} (R)$
over
$\tilde{\cal J}$.
}
\end{thm}

\begin{thm}
{\rm
On
$\tilde{\cal J} \times_{{\rm Spec} {\cal O}_{K}} {\cal C}$
we will consider the hermitian line bundle
$$\overline{\cal F} := (c \times id)^{*} \overline{\cal M} ~~,$$
where
$c: \tilde{\cal J} \longrightarrow {\cal J}$
denotes here the extension of
$C^{(g)} \longrightarrow J$
(the blow-down morphism).
}
\newline
{\bf Facts.}
{\rm a)}
${\cal F}|_{J \times C} = {\cal E} \otimes {\cal O} (R)$.
\newline
{\rm b)}
One has
$\det R \pi_{*} {\cal F} \cong {\cal O}_{\tilde{\cal J}}$.
\newline
{\rm Note here,
$\tilde{\cal J}$
is not regular, so we do not know whether
$\pi: \tilde{\cal J} \times_{{\rm Spec} {\cal O}_{K}} {\cal C} \longrightarrow
\tilde{\cal J}$
has finite Tor-dimension. Thus
$\det R \pi_{*}$
does may be not exist as a functor, but for line bundles, coming by base change
from
${\cal J} \times_{{\rm Spec} {\cal O}_{K}} {\cal C}$,
the definition makes sense.
}
\end{thm}

\begin{thm}
{\bf Remark.}
{\rm
All in all we obtain a decomposition
$$\overline{\cal F}^{\otimes N} \cong \overline{\cal K} \otimes \pi^{*}
\overline{\cal R}$$
of hermitian line bundles, where
${\cal K}$
extends
${\cal E}^{\otimes N}$,
the
$N$-th
power of the tautological line bundle on
$C^{(g)} \times C$.
}
\end{thm}

\begin{thm}
{\bf Remark.}
{\rm Any line bundle of degree
$g$
on
$C$
gives an
${\cal O}_{K}$-valued
point
$y: {\rm Spec} ~ {\cal O}_{K} \longrightarrow {\cal J}$
and
${\rm Spec} ~ {\cal O}_{K} \times_{{\cal J}, y} \tilde{\cal J}$
will be proper over
${\rm Spec} ~ {\cal O}_{K}$.
So at least for some finite field extension
$L/K$
there will be an
${\cal O}_{L}$-valued
point
$\underline{y}: {\rm Spec} ~ {\cal O}_{L} \longrightarrow \tilde{\cal J}$
lifting
$y$.
Proposition \ref{Ara}.c) gives
$$\stackrel{\wedge}{\deg} \left( \det R \pi_{*} \left( {\cal F} |_{\overline{\{
y \} } \times_{{\rm Spec} {\cal O}_{K}} {\cal C}}, h_{{\cal F},y} \right) , \|
. \|_{Q,h} \right) = 0 ~~.$$
}
\end{thm}

\section{Decomposition into two summands}
\begin{thm}
{\rm In this section we will restrict to the case that
${\cal C}$
is {\it semistable}, i. e.
$\pi: {\cal C} \longrightarrow {\rm Spec} ~ {\cal O}_{K}$
is smooth up to codimension
}
$2$.
\newline
{\bf Lemma.}
{\rm a)}
$\tilde{\cal J} \times_{{\rm Spec} {\cal O}_{K}} {\cal C}$
is a normal scheme.
\newline
{\rm b)}
$\tilde{\cal J}$
is quasi-projective over
${\cal O}_{K}$.
\newline
{\bf Proof.}
{\rm a)
$\tilde{\cal J}$
is normal, so
$\tilde{\cal J} \times_{{\rm Spec} {\cal O}_{K}} {\cal C}^{\rm smooth}$
is normal by [SGA1, Expos\'e I, Corollaire 9.10]. In particular
$\tilde{\cal J} \times_{{\rm Spec} {\cal O}_{K}} {\cal C}$
is regular in codimension
$1$.
Further
$\pi: \tilde{\cal J} \times_{{\rm Spec} {\cal O}_{K}} {\cal C} \longrightarrow
\tilde{\cal J}$
is flat with one dimensional fibers. By [EGA IV, Corollaire 6.4.2]
$\tilde{\cal J} \times_{{\rm Spec} {\cal O}_{K}} {\cal C}$
is Cohen-Macaulay in codimension
$2$.
\newline
b)
$\cal J$
is quasi-projective over
${\cal O}_{K}$
by [CS, Chapter VIII by M. Artin, \S 4] and blow-ups are projective morphisms.
}
\begin{center}
$\Box$
\end{center}
\end{thm}

\begin{thm}
{\bf Remark.}
{\rm Let
$\underline{y}: {\rm Spec} ~ {\cal O}_{L} \longrightarrow \tilde{\cal J}$
be an
${\cal O}_{L}$-valued
point lifting an
${\cal O}_{K}$-valued
point \linebreak
$y: {\rm Spec} ~ {\cal O}_{K} \longrightarrow {\cal J}$.
Then
\begin{eqnarray*}
h_{x, \omega} \left( {\cal M}|_{y \times_{{\rm Spec} {\cal O}_{K}} {\cal C}}
\right) & = & \frac{\scriptstyle 1}{\scriptstyle [L:K] N}
\stackrel{\wedge}{\deg} \pi_{*} \left[ \stackrel{\wedge}{c_{1}} \left(
\overline{\cal F}^{\otimes N} |_{\underline{y} \times_{{\rm Spec} {\cal O}_{K}}
{\cal C}} \right) \cdot (x,g_{x}) \right] \\
 & = & \frac{\scriptstyle 1}{\scriptstyle [L:K] N} \stackrel{\wedge}{\deg}
\pi_{*} \left[ \stackrel{\wedge}{c_{1}} \left( \overline{\cal K}
|_{\underline{y} \times_{{\rm Spec} {\cal O}_{K}} {\cal C}} \right) \cdot
(x,g_{x}) \right] \\
 & + & \frac{\scriptstyle 1}{\scriptstyle [L:K] N} \stackrel{\wedge}{\deg}
\pi_{*} \left[ \stackrel{\wedge}{c_{1}} \left( \pi^{*} \overline{\cal R}
|_{\underline{y} \times_{{\rm Spec} {\cal O}_{K}} {\cal C}} \right) \cdot
(x,g_{x}) \right] ~~.
\end{eqnarray*}
We note here, on
${\rm Spec} ~ {\cal O}_{L} \times_{{\rm Spec} ~ {\cal O}_{K}} {\cal C}$
being in general a singular scheme,
$\pi_{*}$
of an intersection with an arithmetic Chern class is defined using an embedding
into a regular scheme, where the line bundle comes from by base change (Remark
\ref{sing}, \cite{Fulton 75}.) The first equality comes from projection
formula. Note that
$x$
means here an
${\cal O}_{L}$-valued
point of
${\rm Spec} ~ {\cal O}_{L} \times_{{\rm Spec} ~ {\cal O}_{K}} {\cal C}$
whose push-forward to
$\cal C$
is
$[L : K] (x)$.
}
\end{thm}

\begin{thm}
{\rm Let us investigate the first summand. We have
${\cal K} |_{C^{(g)} \times C} = {\cal E}^{\otimes N}$
and this line bundle has a canonical section
$s$,
which can be extended over the finite places. Using this section we obtain the
arithmetic cycle
$({\rm div} ~ (s), -\log \| s \|^{2})$
representing
$$\stackrel{\wedge}{c_{1}} (\overline{\cal K}) \in ~ \stackrel{\wedge}{CH^{1}}
\left( \tilde{\cal J} \times_{{\rm Spec} {\cal O}_{K}} {\cal C} \right) ~~.$$
The scheme part
${\rm div} (s)$
of this cycle is an extension of the tautological divisor representing
$c_{1} ({\cal E}^{\otimes N}) \in {\rm CH}^{1} (C^{(g)} \times C)$
(whose restriction to
$\{ (x_{1}, \ldots, x_{g}) \} \times C$
is
$N (x_{1}) + \ldots + N (x_{g})$).
So
${\rm div} ~ (s)$
is the closure of that divisor, possibly plus a finite sum of divisors over the
finite places. We note, that
$\cal K$
is given by that divisor since
$\tilde{\cal J} \times_{{\rm Spec} {\cal O}_{K}} {\cal C}$
is normal. Consequently, if
$\underline{y}$
restricts to the
$L$-valued
point corresponding to the divisor
$D$
on
$C$,
then
$$c_{1} \left( {\cal K} |_{\underline{y} \times_{{\rm Spec} {\cal O}_{K}} {\cal
C}} \right) = \left( \overline{D} \right) + ({\rm correction ~ terms}) ~~,$$
where
$\overline{D}$
denotes the closure of
$N D$
over
$\cal C$
and the correction terms are vertical divisors which (over all the
$y$)
occur only over a finite amount of finite places. Their intersection numbers
with
$(x, g_{x})$,
i. e. with the line bundle
${\cal O} (x)$,
are bounded by
${\rm O} ([L:K])$.

The infinite part
$f$
of
$\stackrel{\wedge}{c_{1}} (\overline{\cal K}) = ({\cal D}, f)$
is a function on
$(C^{(g)} \times C) \backslash {\rm div} ~ (s)$
whose pullback to
$C^{g} \times C$
satisfies all the assumptions of Lemma \ref{Int}. We obtain
\begin{eqnarray*}
\frac{\scriptstyle 1}{\scriptstyle [L:K] N} \stackrel{\wedge}{\deg} \pi_{*}
\left[ \stackrel{\wedge}{c_{1}} \left( \overline{\cal K} |_{\underline{y}
\times_{{\rm Spec} {\cal O}_{K}} {\cal C}} \right) \cdot (x,g_{x}) \right] & =
& \frac{\scriptstyle 1}{\scriptstyle [L:K] N} \left[ \stackrel{\wedge}{\deg}
\left( \stackrel{\wedge}{c_{1}} \left( \overline{{\cal O} (x)} \right) ~ \Big|
{}~ \left( {\cal D} |_{\underline{y} \times_{{\rm Spec} ~ O_{K}} {\cal C}}
\right) \right) ~~~ \ldots \right. \\
 & ~ & ~~~\ldots ~~~ + \left. \frac{\scriptstyle 1}{\scriptstyle 2}
\sum_{\sigma : L \hookrightarrow {\Bbb {C}}} \int_{C ({\Bbb{C}})} f_{D}
\omega_{x} \right] \\
 & = & \frac{\scriptstyle 1}{\scriptstyle N} h_{\overline{\cal S}} (N D) + {\rm
O} (1) \\
 & = & h_{\overline{\cal S}} (D) + {\rm O} (1) ~~.
\end{eqnarray*}
Note, for the first equation we used the symmetry of the intersection form for
hermitian line bundles (Remark \ref{sing}). The denominator
$[L:K]$
disappears by [BoGS, formula (3.1.8)].
}
\end{thm}

\begin{thm}
{\rm The second summand is simpler. One has
\begin{eqnarray*}
\stackrel{\wedge}{c_{1}} \left( \pi^{*} \overline{\cal R} |_{\underline{y}
\times_{{\rm Spec} {\cal O}_{K}} {\cal C}} \right) \cdot (x, g_{x}) & = &
\pi^{*} \stackrel{\wedge}{c_{1}} \left( {\cal R} |_{\underline{y}} \right)
\cdot (x, g_{x}) \\
  & = & \stackrel{\wedge}{c_{1}} \left( {\cal R} |_{\underline{y}} \right) +
\left( 0, g_{x} \omega_{R} (\underline{y}) \right) ~~,
\end{eqnarray*}
when one identifies
$\tilde{\cal J} \times_{{\rm Spec} ~ {\cal O}_{K}} \overline{ \{ x \} }$
with
$\tilde{\cal J}$.
The integral
$\int_{C ({\Bbb C})} ~ g_{x} \omega_{R} ( . )$
depends smoothly on the parameter,
in particular it is bounded. So the push-forward of the right summand is
bounded by
${\rm O} ([L:K])$.
On the other hand
$\pi_{*} \stackrel{\wedge}{c_{1}} \left( {\cal R} |_{\underline{y}} \right) =
\left( \stackrel{\wedge}{c_{1}} (\overline{\cal R}) ~ \Big| ~ \underline{y}
\right) $,
where the last term is defined by embedding
$\tilde{\cal J}$
into a scheme
$P$
smooth and projective over
${\cal O}_{K}$
\cite{Fulton 75}. Note here we use
$\tilde{\cal J}$
is quasi-projective. Thus Lemma \ref{class} gives
$$\stackrel{\wedge}{\deg} \left( \stackrel{\wedge}{c_{1}} (\overline{\cal R}) ~
\Big| ~ \underline{y} \right) = \stackrel{\wedge}{\deg} \left(
\stackrel{\wedge}{c_{1}} (\overline{{\cal R}_{P}}) ~ \Big| ~ \iota_{*}
(\underline{y}) \right) = h_{\overline{{\cal R}_{P}}} \left( \iota_{*}
(\underline{y}) \right) = h_{R_{P}} \left( \iota_{*} D \right) = h_{R} (D)
{}~~,$$
where
$D$
is the divisor corresponding to the restriction of
$\underline{y}$
to
$C^{(g)}$.
}
\end{thm}

\begin{thm}
\label{stab}
{\rm We obtain}
\newline
{\bf Proposition.}
Assume
$\cal C$
is semistable and
${\cal L} = {\cal O} (\overline{D})$,
where
$\overline{D}$
is the closure of some divisor on
$C$.
Then Theorem \ref{main} is true.
\newline
{\bf Proof.}
{\rm By Corollary \ref{two} this is now proven for line bundles coming by
restriction from
$\cal M$.
This way one can realize the line bundles
${\cal O} (D)$
on the generic fiber
$C$
for arbitrary divisors
$D$
(defined over
$K$)
of degree
$g$
over
$C$.
Consider the degrees
$$\deg {\cal M} |_{y \times_{{\rm Spec} {\cal O}_{K}} {\cal C}_{{\goth p},
i}}$$
for
${\cal O}_{K}$-valued
points
$y$
of
$\cal J$,
where
${\cal C}_{{\goth p}, i}$
denote the irreducible components of the special fiber
${\cal C}_{\goth p}$.
They are even defined for
$\overline{{\cal O} / {\goth p}}$-valued
points, where the bar denotes algebraic closure here, and are locally constant
over the special fiber
${\cal J}_{\goth p}$.
In particular they are bounded since the N\'eron model of an abelian variety is
of finite type. The Proposition follows from Lemma \ref{degr}.

}
\begin{center}
$\Box$
\end{center}
\end{thm}

\section{End of the proof}
\begin{thm}
\label{blow}
{\bf Lemma.}
{\rm Let
${\cal C} / {\cal O}_{K}$
be a regular projective arithmetic surface and
$p: \tilde{\cal C} \longrightarrow {\cal C}$
be a blow-up of one point. Then
$$h_{x, \omega} ({\cal L}) = h_{x, \omega} (p^{*} {\cal L}) ~~,$$
where
$x \in {\cal C} ({\cal O}_{K}) = \tilde{\cal C} ({\cal O}_{K})$
and
$\cal L$
is a line bundle with
$\chi ({\cal L}) \neq 0$.
\newline
{\bf Proof.}
Obviously
$p_{*} p^{*} {\cal L} = {\cal L}$
and [SGA6, Expos\'e VII, Lemma 3.5] gives
$R^{i} p_{*} p^{*} {\cal L} = 0$
for
$i \geq 1.$
In particular
$R p_{*} (p^{*} {\cal L}) = {\cal L}$,
$R (\pi p)_{*} (p^{*} {\cal L}) = R \pi_{*} {\cal L}$
and
$\det R (\pi p)_{*} (p^{*} {\cal L}) = \det R \pi_{*} {\cal L}$.
\linebreak This means that
$\cal L$
and
$p^{*} {\cal L}$
get identical distinguished metrics and therefore \linebreak
$\stackrel{\wedge}{c_{1}} \left( p^{*} {\cal L}, \| . \|_{p^{*} {\cal L}}
\right) = p^{*} \stackrel{\wedge}{c_{1}} \left( {\cal L}, \| . \|_{{\cal L}}
\right) $.
On the other hand
$p^{*} (x, g_{x}) = (x, g_{x}) + ({\rm exceptional ~ divisor})$,
but an exceptional divisor intersects trivially with cycles coming from
downstairs. Consequently,
\begin{eqnarray*}
h_{x, \omega} (p^{*} {\cal L}) & = & \stackrel{\wedge}{\deg} (\pi p)_{*} \left[
p^{*} \stackrel{\wedge}{c_{1}} ({\cal L}, \| . \|_{\cal L}) \cdot p^{*} (x,
g_{x}) \right] \\
 & = & \stackrel{\wedge}{\deg} (\pi p)_{*} p^{*} \left[
\stackrel{\wedge}{c_{1}} ({\cal L}, \| . \|_{\cal L}) \cdot (x, g_{x}) \right]
\\
 & = & \stackrel{\wedge}{\deg} \pi_{*} \left[ \stackrel{\wedge}{c_{1}} ({\cal
L}, \| . \|_{\cal L}) \cdot (x, g_{x}) \right] ~~~~~~~~~~~~~~~~~~~~{\rm
"projection ~ formula"}\\
 & = & h_{x, \omega} ({\cal L}) ~~.
\end{eqnarray*}
}
\begin{center}
$\Box$
\end{center}
\end{thm}

\begin{thm}
{\bf Corollary.}
\label{model}
{\rm (Change of model.)}
\newline
Let
${\cal C}_{1}, {\cal C}_{2} / {\cal O}_{K}$
be two regular projective models of the curve
$C/K$
of genus
$g$.
Then, for divisors
$D$
of degree
$g$
on
$C$
$$h_{x, \omega} \left( {\cal O}_{{\cal C}_{1}} (\overline{D}) \right) = h_{x,
\omega} \left( {\cal O}_{{\cal C}_{2}} (\overline{D}) \right) + {\rm O} (1)
{}~~,$$
where
$\overline{D}$
denotes the closure of
$D$
in
${\cal C}_{1}$,
respectively
${\cal C}_{2}.$
\newline
{\bf Proof.}
{\rm By [Li, Theorem II.1.15] one is reduced to the case of the blow-up of one
point
$p: {\cal C}_{2} \longrightarrow {\cal C}_{1}$.
By Lemma \ref{blow} we have to bound the difference
$h_{x, \omega} \left( {\cal O}_{{\cal C}_{2}} (\overline{D}) \right) - h_{x,
\omega} \left( p^{*} {\cal O}_{{\cal C}_{1}} (\overline{D}) \right)$.
$\overline{D}$
will meet the point blown up
$i$
times
($0 \leq i \leq g$).
We get an exact sequence
$$0 \longrightarrow {\cal O}_{{\cal C}_{2}} (\overline{D}) \longrightarrow
p^{*} {\cal O}_{{\cal C}_{1}} (\overline{D}) \longrightarrow {\cal O}_{E^{i}}
\longrightarrow 0 ~~,$$
where
$E$
is the exceptional curve and
$E^{i}$
denotes its
$i$-th
infinitesimal neighbourhood. But now the assertion is a direct consequence of
Lemma \ref{red}.
}
\begin{center}
$\Box$
\end{center}
\end{thm}

\begin{thm}
{\bf Lemma.}
{\rm (Change of base field.)}
\newline
Let
${\cal C} / {\cal O}_{K}$
be a regular arithmetic surface, generically of genus
$g$,
$L/K$
a finite field extension and
$$p: {\cal C}^{'} = \overline{{\cal C} \times_{{\rm Spec} ~ {\cal O}_{K}} {\rm
Spec} ~ {\cal O}_{L}} \longrightarrow {\cal C}$$
be some resolution of singularities of the base change to
${\cal O}_{L}$.
Then, for divisors
$D$
of degree
$g$
on
$C = {\cal C} \times_{{\rm Spec} ~ {\cal O}_{K}} {\rm Spec} ~ K$,
$$h_{x, \omega} \left( {\cal O}_{{\cal C}^{'}} (\overline{p^{*} D}) \right) =
[L:K] \cdot h_{x, \omega} \left( {\cal O}_{\cal C} (\overline{D}) \right) +
{\rm O} (1) ~~.$$
\newline
{\bf Proof.}
{\rm
$p$
is a composition of blow-ups and finite morphisms [CS, Chapter XI by M. Artin].
\linebreak Using the first formulas in the proof of Lemma \ref{blow}
successively we obtain
\linebreak
$R p_{*} p^{*} {\cal O} (\overline{D}) = {\cal O} (\overline{D})$
and
$\det R (\pi p)_{*} p^{*} \left( {\cal O} (\overline{D}) \right) = \det R
\pi_{*} {\cal O} (\overline{D})$
such that
${\cal O} (\overline{D})$
and
$p^{*} {\cal O} (\overline{D})$
get identical distinguished metrics. Here it follows
$$h_{x, \omega} \left( p^{*} {\cal O} (\overline{D}) \right) = [L:K] \cdot
h_{x, \omega} \left( {\cal O} (\overline{D}) \right) ~~,$$
since
$p$
is a morphism of degree
$[L:K]$
and the projection formula gives
$p_{*} p^{*} Z = [L:K] Z$.
$p^{*} {\cal O} (\overline{D})$
and
${\cal O} (\overline{p^{*} D})$
differ by a limited combination of the exceptional divisors such that the
assertion follows from Lemma \ref{red}.
}
\begin{center}
$\Box$
\end{center}
\end{thm}

\begin{thm}
\label{hori}
{\bf Proposition.}
For line bundles
${\cal L} = {\cal O} (\overline{D})$,
where
$\overline{D}$
is the closure of some divisor of degree
$g$
on
$C$,
Theorem \ref{main} is true.
\newline
{\bf Proof.}
{\rm By [AW, Corollary 2.10] there is a stable model for
$C \times_{{\rm Spec} ~ K} {\rm Spec} ~ L$
after some finite field extension
$L/K$.
The assertion follows from Proposition \ref{stab}.
}
\begin{center}
$\Box$
\end{center}
\end{thm}

\begin{thm}
{\bf Proposition.}
Theorem \ref{main} is true.
\newline
{\bf Proof.}
{\rm This is a direct consequence of Proposition \ref{hori} and Lemma
\ref{degr}.
}
\begin{center}
$\Box$
\end{center}
\end{thm}

\section{Some technical Lemmata}
\begin{thm}
\label{fibe}
{\bf Lemma.}
{\rm (Fibers do not change the height.)}
\newline
If
$\cal L$
is a line bundle on
${\cal C} / {\cal O}_{K}$
with
$\chi ({\cal L}) \neq 0$,
then
$$h_{x, \omega } \left( {\cal L} \otimes {\cal O} (\goth{p}) \right) = h_{x,
\omega } ({\cal L})$$
for every prime ideal
$\goth{p} \subseteq {\cal O}_{K}$.
\newline
{\bf Proof.}
{\rm One has
${\cal O} (\goth{p}) = \pi^{*} (\goth{p}^{-1})$,
hence by projection formula
$$\det R \pi_{*} \left( {\cal L} \otimes {\cal O} (\goth{p}) \right) \cong \det
R \pi_{*} {\cal L} \otimes {\cal O} ({\goth p})^{- \chi ({\cal L})} ~~.$$
Let
$\| . \|$
be one of the distinguished metrics on the line bundle
${\cal L}_{\Bbb{C}}$
on
$\coprod_{\sigma: K \hookrightarrow \Bbb{C}} C(\Bbb{C})$.
\linebreak
We put
$\| . \|_\goth{p} = C \cdot \| . \|$
for a distinguished hermitian metric on
$({\cal L} \otimes {\cal O} (\goth{p}))_{\Bbb{C}} = {\cal L}_{\Bbb{C}}$.
It follows
\linebreak
$h_{Q, \det R \pi_{*} ({\cal L} \otimes {\cal O} (\goth{p}))} = C^{\chi ({\cal
L})} \cdot h_{Q, \det R \pi_{*} {\cal L}}$
and
\begin{eqnarray*}
\stackrel{\wedge}{\deg} \left( \det R \pi_{*} ({\cal L} \otimes {\cal O}
(\goth{p})), h_{Q, \det R \pi_{*} ({\cal L} \otimes {\cal O} (\goth{p}))}
\right) & = & \stackrel{\wedge}{\deg} \Big( \det R \pi_{*} {\cal L}, h_{Q, \det
R \pi_{*} {\cal L}} \Big) \\  & + &  \chi ({\cal L}) \Big[ [K : {\Bbb{Q}}] \log
C - \log (\sharp {\cal O} / \goth{p}) \Big] ~~.
\end{eqnarray*}
Thus a distinguished hermitian metric on
$({\cal L} \otimes {\cal O} (\goth{p}))_{\Bbb{C}}$
can be given by
$\| . \|_\goth{p} = ( \sharp {\cal O} / \goth{p})^{\frac{1}{[K : {\Bbb Q}]}}
\cdot \| . \|$
and it follows
$$\stackrel{\wedge}{c_{1}} ({\cal L} \otimes {\cal O} (\goth{p}), \| .
\|_\goth{p}) = ~ \stackrel{\wedge}{c_{1}} ({\cal L}, \| . \|) + \pi^{*} \left(
\goth{p} ; - \frac{\scriptstyle 2}{\scriptstyle [K : \Bbb{Q}]} \log (\sharp
{\cal O} / \goth{p}) ,\ldots ,- \frac{\scriptstyle 2}{\scriptstyle [K :
\Bbb{Q}]} \log (\sharp {\cal O} / \goth{p}) \right) ~~.$$
But the arithmetic cycle
$\left( \goth{p} ; - \frac{2}{[K : \Bbb{Q}]} \log (\sharp {\cal O} / \goth{p}),
\ldots, - \frac{2}{[K : \Bbb{Q}]} \log (\sharp {\cal O} / \goth{p}) \right) \in
{}~ \stackrel{\wedge}{{\rm CH}^{1}}({\rm Spec} ~ {\cal O}_{K})$
vanishes after multiplication with the class number
$\sharp {\rm Pic} ~ ({\rm Spec} ~ {\cal O}_{K})$,
hence it is torsion and therefore numerically trivial.
}
\begin{center}
$\Box$
\end{center}
\end{thm}

\begin{thm}
\label{red}
{\bf Lemma.}
Let
$F$
be some vertical divisor on
${\cal C} / {\cal O}_{K}$.
Then, for line bundles
${\cal L} / {\cal C}$,
fiber-by-fiber of degree
$g$,
$$h_{x, \omega} ({\cal L} (F)) = h_{x, \omega} ({\cal L}) + {\rm O} (1) ~~.$$
{\bf Proof.}
{\rm By Lemma \ref{fibe} we may assume that
$E := -F$
is effective. Using induction we are reduced to the case
$E$
is an irreducible curve. We have a short exact sequence
$$0 \longrightarrow {\cal L} (F) \longrightarrow {\cal L} \longrightarrow {\cal
L}_{E} \longrightarrow 0$$
inducing the isomorphism
$$\det R \pi_{*} {\cal L} (F) \cong \det R \pi_{*} {\cal L} \otimes (\det R
\pi_{*} {\cal L}_{E})^{\vee} ~~.$$
But
$\det R \pi_{*} {\cal L}_{E}$
depends only on the Euler characteristic of
${\cal L}_{E}$
and for the degree of that bundle there are only
$g+1$
possibilities. So up to numerical equivalence there are only
$g+1$
possibilities for
$$\stackrel{\wedge}{c_{1}} \left( {\cal L} (F), \| . \|_{{\cal L} (F)} \right)
- \stackrel{\wedge}{c_{1}} \Big( {\cal L}, \| . \|_{\cal L} \Big) ~~,$$
where
$\| . \|_{\cal L}$
and
$\| . \|_{{\cal L} (F)}$
denote distinguished hermitian metrics.
}
\begin{center}
$\Box$
\end{center}
\end{thm}

\begin{thm}
\label{degr}
{\bf Lemma.}
Consider line bundles
${\cal L}$,
generically of degree
$g$
on
${\cal C}$,
equipped with a section
$s \in \Gamma (C, {\cal L}_{C})$
over the generic fiber, and assume the degrees
$\deg {\cal L} |_{{\cal C}_{{\goth p}, i}}$
of the restrictions of
$\cal L$
to the irreducible components of the special fibers to be fixed. Then
$$h_{x, \omega} ({\cal L}) = h_{x, \omega} \left( {\cal O} \left(
\overline{{\rm div} (s)} \right) \right) + {\rm O} (1) ~~.$$
{\bf Proof.}
{\rm We have
${\cal L} = {\cal L}^{'} (E)$,
where
${\cal L}^{'} = {\cal O} \left( \overline{{\rm div} (s)} \right)$
is a line bundle induced by a horizontal divisor and
$E$
is a vertical divisor. By Lemma \ref{fibe} we may assume
$E$
to be concentrated in the reducible fibers of
$\cal C$.
So, using induction, let
$E$
be in one such fiber
${\cal C}_{{\goth p}}$.
Then for the degrees
$\deg {\cal O} (E) |_{{\cal C}_{{\goth p}, i}}$
there are only finitely many possibilities. But by [Fa, Theorem 4.a)] the
intersection form on
${\cal C}_{\goth p}$
\pagebreak
is negative semi-definite where only multiples of the fiber have square
$0$. \linebreak
Hence, for
$E$
there are only finitely many possibilities up to addition of the whole fiber,
which does not change the height. Lemma \ref{red} gives the claim.
}
\begin{center}
$\Box$
\end{center}
\end{thm}

\begin{thm}
{\bf Lemma.}
\label{Int}
Let
$X$
be a compact Riemann surface and
$g \in \Bbb{N}$
be a natural number. Denote by
$\Delta$
the diagonal in
$X \times X$,
by
$\delta_{M}$
the
$\delta$-distribution
defined by
$M$
and by
$\pi_{i}: X^{g} \times X \longrightarrow X$
(resp.
$\pi_{i,g+1}: X^{g} \times X \longrightarrow X \times X$)
the canonical projection on the
$i$-th
component (resp. to the product of the
$i$-th
and
$(g+1)$-th
component.) Further let
$$f: (X^{g} \times X) \backslash \bigcup_{i=1}^{g} \pi_{i,g+1}^{-1} (\Delta)
\longrightarrow {\Bbb{C}}$$
be a smooth function such that the restriction of
$$-d_{X} d_{X}^{c} f + \delta_{\Delta} \circ \pi_{1,g+1} + \ldots +
\delta_{\Delta} \circ \pi_{g,g+1} = \rho ~~,$$
to
$\{ (x_{1}, \ldots, x_{g}) \} \times X$
is a smooth
$(1,1)$-form
smoothly varying with
$(x_{1}, \ldots, x_{g})$.
Let
$\omega$
be a smooth
$(1,1)$-form
on
$X$.
Then
$$\int_{X} f(x_{1}, \ldots ,x_{g}, \cdot) \omega$$
depends smoothly on
$(x_{1}, \ldots, x_{g}) \in X^{g}.$
\newline
{\bf Proof.}
{\rm Without restriction we may assume
$\int_{X} \omega = 1$.
Then, for any
$x \in X$
there exists a function
$h \in C^{\infty} (X \backslash \{ x \})$,
having a logarithmic singularity in
$x$,
such that
$\omega = -dd^{c} h + \delta_{x}$.
It follows
\begin{eqnarray*}
\int_{X} f(x_{1}, \ldots, x_{g}, \cdot) \omega & = & -\int_{X} f(x_{1}, \ldots,
x_{g}, \cdot) dd^{c} h + f(x_{1}, \ldots, x_{g}, x) \\
 & = & -\int_{X} \Big( d_{X} d_{X}^{c} f(x_{1}, \ldots ,x_{g}, \cdot) \Big) h +
f(x_{1}, \ldots, x_{g}, x) \\
 & = & \int_{X} \rho (x_{1}, \ldots, x_{g}, \cdot) h - h(x_{1}) - \ldots -
h(x_{g}) + f(x_{1}, \ldots, x_{g}, x) \\
 & = & \int_{X} \rho (x_{1}, \ldots, x_{g}, \cdot) h - \Big[ h(x_{1}) - G(x,
x_{1}) \Big] - \ldots - \Big[ h(x_{g}) - G(x, x_{g}) \Big] \\
 & - & \Big[ G(x, x_{1}) + \ldots + G(x, x_{g}) - f(x_{1}, \ldots, x_{g}, x)
\Big] ~~,
\end{eqnarray*}
where
$G$
is the Green's function of
$X$.
Because
$h$
has only a logarithmic singularity it is allowed to differentiate under the
integral sign. So the integral is smooth. The other summands are solutions of
equations of the form
$dd^{c} F = \sigma$
with a smooth
$(1,1)$-form
$\sigma$
on
$X$
satisfying
$\int_{X} \sigma = 0$
(in
$x_{1}, \ldots, x_{g}$,
respectively
$x$).
Since
$dd^{c}$
is elliptic, these solutions exist as smooth functions and are unique up to
constants. In particular, also the last summand must depend smoothly on
$(x_{1}, \ldots, x_{g})$,
even when some of the
$x_{i}$
equal
$x$.
Note that the symmetry of the Green's function is used here essentially.
\begin{center}
$\Box$
\end{center}
}
\end{thm}
{}~
\newline
{\footnotesize {\bf Acknowledgement.} When doing this work, the author had
fruitful discussions with U. Bunke (Berlin), who explained him much of the
analytic part of the theory. He thanks him warmly.}

\newpage
{}~
\vspace*{-1.20truecm}
\thispagestyle{myheadings}
\markright{\rm ~~~~~~~~~~~~~~~~~~~~~~~~~~~~~~~~~~~~~~~~~~~~~~~~~~~~~~~~~Jahnel}

\thispagestyle{myheadings}
\markright{\rm ~~~~~~~~~~~~~~~~~~~~~~~~~~~~~~~~~~~~~~~~~~~~~~~~~~~~~~~~~Jahnel}
\end{document}